\documentclass[12pt]{article}

\usepackage{amsmath}
\usepackage{amssymb}
\usepackage{amsthm,bm}
\usepackage{amsfonts}
\usepackage{ulem}
\usepackage{cancel}
\usepackage{graphicx}
\usepackage{dsfont} 
\usepackage{soul}
\usepackage[mathscr]{eucal}

\usepackage{txfonts}

\usepackage{bbm,epsfig,graphics,epic,color,rotating,xcolor, eucal, eufrak}
\numberwithin{equation}{section}

\newcommand{\R}{{\mathbb R}}

\newcommand{\A}[1]{{\mathfrak #1}}
\newcommand{\Ab}[1]{\pmb{\mathfrak #1}}

\newcommand{\ga}{\gamma}

\newcommand{\Proof}{{\textit {Proof} }}

\newcommand{\wt}[1]{\widetilde{#1}}
\newcommand{\wh}[1]{\widehat{#1}}

\newcommand{\ed}{\mathrm{d}}

\newcommand{\reff}[1]{(\ref{#1})}

\renewcommand{\k}{\varkappa}
\newcommand\e{\varepsilon}
\newcommand\la{\lambda}
\newcommand{\1}{{\mathbbm{1}}}
\renewcommand\phi{\varphi}

\newcommand{\be}{\begin{equation}}
\newcommand{\ee}{\end{equation}}
\newcommand{\bel}[1]{\begin{equation}\label{#1}}
\newcommand{\bea}{\begin{eqnarray}}
\newcommand{\eea}{\end{eqnarray}}
\newcommand{\balign}{\begin{aligned}}
\newcommand{\ealign}{\end{aligned}}
\newcommand{\ba}{\begin{array}}
\newcommand{\ea}{\end{array}}
\newcommand{\bfig}{\begin{figure}}
\newcommand{\efig}{\end{figure}}

\newcommand{\PP}{{\mathbf P}}

\newtheorem{theorem}{Theorem}[section]
\newtheorem{lemma}[theorem]{Lemma}
\newtheorem{proposition}[theorem]{Proposition}

\newtheorem{definition}[theorem]{Definition}
\newtheorem{remark}[theorem]{Remark}

\theoremstyle{definition}

\title{Sojourn times of Markov symmetric processes in continuous time.}

\author{Pechersky E.A.$^1$, Presman E.L.$^2$, Yambartsev A.A.$^3$
\footnote{AY thanks FAPESP for the finantial support via the grant 2017/10555-0.}
}
\date{}

\begin{document}

\maketitle {\footnotesize
	
	\noindent $^1$Dobrushin´s Laboratory, Institute of Information Transmission Problems of Russian Academy of Science, Moscow, Russia.
	 \\ E-mail: pech@iitp.ru
	
	\noindent $^2$Central Economics and Mathematics Institue of Russian Academy of Sciences (CEMI RAS), Moscow, Russia.
	 \\ E-mail: presman@cemi.rssi.ru
	
	\noindent $^3$Department of Statistics, Institute of Mathematics
	and Statistics, University of S\~ao Paulo, S\~ao Paulo, Brazil.
	\\ E-mail: yambar@ime.usp.br}

\vspace{0.5cm}


\vspace{1cm}

\hspace{1cm} \textit{In memory of Nikita D. Vvedenskaya}

\begin{abstract}
	The symmetric birth and death process in the integers $\{1, \ldots, N \}$ with linear rates is studied. The process moves slowly and spends more time in the neighborhood of the state 1. It represents our attempt at explaining the asymmetry between amounts of matter and antimatter by the inhomogeneity of the sojourn time. The state of the process reflects a relative frequency of an antimatter amount. The observed matter-antimatter imbalance in the Universe we consider a result of stochastic competition between them. The amount of matter significantly exceeds the amount of antimatter, which corresponds to the lower states of the process, and the path toward matter-antimatter equilibrium can be very long.  
\end{abstract}

\section{Introduction}
 \subsection{Markov process} 
 
 
 \noindent
 \textit{Markov Process.}
 We study a continuous time Markov process $X^N(t)$ with the state space $\mathcal N=\{1,...,N\}$ and piecewise constant paths. We consider two kinds of the initial state, it is either a point from $\mathcal N$ or the stationary distribution \reff{stat.dist} on $\mathcal N$.
 The jumps of the process with non-zero transition rates are $k\to k+1$ and $k\to k-1$ with the rates $\la_{k, k+1}, \la_{k,k-1}$ if $1<k< N$, and $1\to2, N\to N-1$ with the rates $\lambda_{1, 2} = 1, \lambda_{N, N-1}=N$. The dynamics of the process is described by the \textit{site-symmetry} generator with rates $\la_{k,k+1}=\la_{k,k-1}=k$ for $k\neq1,N$, thus, in this case we have the \textit{edge-asymmetric}, $\la_{k,k+1}\ne\la_{k+1,k}$, see \reff{1.2}. It means that the considered processes $X^N$ are space non-homogeneous.  
 
\vspace{0.5cm}
\noindent
\textit{Invariant measures.}
 Note that the corresponding {\it embedded} discrete-time Markov chain $\wh{X}^N$ is a symmetric simple random walk over $\mathcal{N}$ with reflections at the boundaries $1$ and $N$ (see Section \ref{s3.1}). However, the continuous-time process $X^N$ is far from any symmetry. As will be shown, despite the property of symmetry $\la_{k,k+1}=\la_{k,k-1}$ the process $X^N$ spends most of its time in the vicinity of the point $0$. That is, the distribution of $X^N(t)$ is concentrated in the neighborhood of $0$, while the embedded process $\wh{X}^N$ has a uniform distribution over the set $\mathcal N \setminus \{1, N\}$. Both $X^N$ and $\wh{X}^N$ satisfy detailed balance with respect to their stationary distribution.
 
\vspace{0.5cm}
 \noindent
\textit{Law of Large Numbers.}
 We study a sequence $\xi^N$ of the Markov processes obtained by the scaling by $N$ of the original processes $X^N$, that is $\xi^N:=\frac{X^N}N \in \A{N}^N=\{\frac1N,\frac2N,...,1\}$. 
 We show that if the initial state is one-pointed then the mean value of $\xi^N$ tends to  the constant  path equal to the limiting initial value $\xi^N(0)$ with growing $N$.
 And if the initial state is sampled according to the stationary distribution then the limit of the stationary distribution of $\xi^N$ shows the effect of the concentration near the path identically equal to 0. 
 
\vspace{0.5cm}
 \noindent
\textit{Large Deviations.}
 Finally, the rate function of the large deviations is obtained by the method described in \cite{Kurtz2}. For any given $x_0, x_T \in [0,1]$ we find the optimal path $\gamma(t), t\in [0,T]$ starting at $\gamma(0)=x_0$ and taking value $x_T$ at the end of time interval $\gamma(T)=x_T$. We show that the optimal path is a parabola.
 
 
 \subsection{Motivations} 

 


The vast majority of laws in physics are related to a certain symmetry. Often, symmetries imply the existence of conservation laws, such as the well-known conservation laws in mechanics and the conservation of the electric charge. 

However, the most impressive discoveries are associated with asymmetric manifestations in nature. Many asymmetries have been found in the microcosm that obeys quantum laws.  
 
One well-known phenomenon is the asymmetry between matter and antimatter or the matter-antimatter imbalance. It is an impressive problem that has not yet been solved. The attempt to find an answer to this problem is based on the search for a violation of some symmetry in physical law, such as the violation of the law of conservation of baryon charges (see, for example, \cite{GorRub}, \cite{Rub}, \cite{Sa}, \cite{TK}). According to \cite{Sath} this asymmetry is a reason for the matter-antimatter imbalance.
 
There is another reason for the asymmetry: randomness inherent in the system. Randomness can cause fluctuations, often quite large.
 
Very often, the asymmetry of physical processes caused by the random behavior of the system expresses itself through a phenomenon of phase transition. The phase transition is associated with very different behaviors of the same physical system on opposite sides of a ``critical point'' of some system parameter. Usually, a phase transition is related to the randomness created by fluctuations in the system or the system environment.

Here we propose another reason for the asymmetry. It is based in fluctuations modeled as a spatially inhomogeneous continuous-time Markov process. The asymmetry arises due to distinctions in the occupation time spent by the Markov process in different states.
 
The state-space of the process under consideration is a finite set $\mathcal N=\{1,...,N\}$, and the process dynamics creates piecewise constant trajectories with jumps $+1$ or $-1$. The transition probabilities to the left and right neighbors from any state $2<i<N-1$ are equal, but the transition rate linearly increases with the growth of $i$. The process exhibits the asymmetry related to the sojourn time. In the area close to $1\in\mathcal N$ the process spends much more time than in the states near the $N$. Assuming the process describes the particle dynamics, it means that if a particle enters a region close to $1\in\mathcal N$ due to the fluctuations, it takes a long time before leaving this region.
 
This work is the first step toward explaining the asymmetry between matter and antimatter due to the fluctuations created by a stochastic symmetric competition between matter and antimatter. We assume that violations of baryon-antibaryon symmetry arise in the course of Markovian dynamics. The breaking of symmetry occurs when the process path is constant during some time. The longer the constant period, the higher the probability of violation. 

This effect can be sharper if we consider Markov dynamics out of the average in the corresponding direction. Such an event may have a small probability. However, it may strengthen the studied effect. This study belongs to the area of theory probability called \textit{large deviations theory}.
   
Another motivation comes from our study of large deviations for birth-death processes with polynomial transition rates. In the series of papers, \cite{MPY}, \cite{VLSY1}, \cite{VLSY2}, \cite{VLSY3}, we consider the birth-death processes where the transition rates $\lambda_{i,i+1}$ and $\lambda_{i,i-1}$ depends on the state $i$ polynomially: $\lambda_{i,i+1}=Pi^{\alpha}$ and $\lambda_{i,i+1}=Qi^{\beta}$. The papers \cite{VLSY1}, \cite{VLSY2}, \cite{VLSY3} explore the method of exchange of measures suggested in \cite{MPY} in the cases when $\alpha \ne \beta$ or $\alpha = \beta, P\ne Q$. But all our attempts have failed to extend the method for the case $\alpha = \beta, P=Q$. 

The model considered here corresponds to the birth-death processes (on the interval $0, 1, \ldots, N$) with $\alpha = \beta=1, P=Q$, where the method of \cite{Kurtz2} was employed.

\section{Model. Generator }\label{1.3} The main subject of our research is the sequence (in $N$, as $N\to\infty$) of Markov processes $(X^N(t), t \in [0,T])$ for some fixed $T>0$. The upper index means that the states of the process belong to the set $\mathcal N=\{1,\ldots, N\}$. The stochastic dynamics of $X^N(t)$ is controlled by the generators $G^N$, which is an operator on the set of probe functions $\mathcal F^N=\{f:\:\mathcal N\to \R\}$:
 \bel{3.1}
G^Nf(m)=\la m \bigl( f(m+1)-f(m) \bigr) \1_{m<N}(m)+\la m \bigl( f(m-1)-f(m)\bigr) \1_{m>1}(m),
\ee
where $\1_{\mathcal B}(\cdot)$ is an indicator for subset $\mathcal B\subset\mathcal N$
   $$\1 _{\mathcal B}(m)=\begin{cases}1,\mbox{if }m\in\mathcal B,\\
0,\mbox{otherwise},
\end{cases}
$$
and the rates $\la_{m,m+1}$ and $\la_{m,m-1}$ are 
\bel{1.2}
\begin{aligned}&
\lambda_{m, m+1} = \lambda_{m, m-1} = \la m, \ m=2, \ldots, N-1, \\
&\lambda_{1,2} = \lambda,\  \lambda_{N,N-1} =  {\lambda}{N}. 
\end{aligned}
\ee

Note that the endpoints of $\mathcal N$ are reflective. The initial state $X^N(0)$ has  some distribution on $\mathcal N$. In particular, it can be a measure concentrated at the one point state $m\in\mathcal N$.

\begin{remark}It can be seen from the form of $G^N$ that the Markov process has jump-wise paths. The jumps are equal to $\pm1$. The intensities of jumps from the state $m$ depend on the value of $m$, and the intensities to the left and right are equal if $m\in\{2,..., N-1\}$.
Thus, we have symmetrical jumps except for the endpoints of $\mathcal N$.
\end{remark}
 
 \subsection{Stationary distribution}
Detailed balance allows finding stationary distributions of the processes. The stationary distribution concentrates near 0 as $N$ increases.

Denote by ${\pi^N(m), m\in\mathcal N}$,  the stationary measure of $X^N(t)$. 
\begin{theorem}\label{22.03.23-1}
\label{4.0}
$X^N(t)$ satisfies  the detailed balance with respect to the stationary distribution. $\pi^N$ is given by 
\begin{equation}\label{stat.dist}
\pi^N(m)=\frac{\frac1m}{\sum_{k=1}^N\frac1k}, \ \ m\in\mathcal N.
\end{equation}
Let ${\mathcal M}=\{1,...,M\}\subseteq\mathcal N$. Let $h(N)$ be a sequence of integers such that 
$$\frac{\ln(h(N))}{\ln(N)}\to1, \mbox{ and }\limsup_{N\to\infty}\frac N{h(N)}>1.$$ 
Assume that $M=h(N)$. Then 
\[
\pi^N(\mathcal M)\to1, \mbox{ as }N\to\infty,
\]
\end{theorem}

\proof The proof uses the detailed balance. The balance equations have the following shape (see \reff{3.1})
\[
\pi^N(m)m=\pi^N(m+1)(m+1),\mbox{ for all  } 1 \le m \le N-1.
\]
For $m>1$ we obtained  
\[
\pi^N(m)=\pi^N(1)\frac1m.
\]
Since $\sum_{m=1}^N\pi^N(m)=1$ we have
\[
\pi^N(1)=\frac1{\sum_{m=1}^N\frac1m},
\]
and
\bel{3.2}
\pi^N(m)=\frac{\frac1m}{\sum_{m=1}^N\frac1m},\ \ m\in \mathcal{N}.
\ee
The probability of the set $\mathcal M=\{1,...,M\}$ is
\bel{3.0}
\pi^N(\mathcal M)=\sum_{m=1}^M\pi^N(m)=\frac{\sum_{m=1}^M\frac1m}{\sum_{m=1}^N\frac1m}.
\ee
Here we use the Euler representation of the partial sum of the Harmonic series $H_K=\sum_{m=1}^K\frac1m=\ln K+\gamma+\e_K$, where $\gamma$ is the Euler–Mascheroni constant and $K\e_K\to\frac12$ as $K\to\infty$.  Then
\bel{4.1} 
\pi^N(\mathcal M)=\frac{\ln M+\gamma+\e_M}{\ln N+\gamma+\e_N}=\frac{\ln h(N)+\gamma+\e_M}{\ln N+\gamma+\e_N}\to 1
\ee
as $N\to\infty$ \hfill$\blacksquare$

\subsection{Embedded chain}\label{s3.1}The stationary distribution of the embedded Markov chain $\wh X$ is 
fundamentally different from the stationary distribution of the  process $X^N$.

Consider the embedded chain $\wh{X}^N(k)$ corresponding to $X^N(t)$ with initial value $\wh{X}^N(0)=m$. The chain $\wh{X}^N(k)$ instantly changes its state immediately upon the occurence of jumps $X^N(t)$, regardless of the time the process spends in the state before the jump. The discrete-time  dynamics is determined by the following probability transition matrix
\[
\begin{pmatrix}
0&\displaystyle{\frac{\la_{1,2}} {R_1}}&0&0&0&...&0&0\\
\displaystyle\frac{\la_{2,1}}{R_2}&0&\displaystyle\frac{\la_{2,3}}{R_2}&0&0&...&0&0\\
0&\displaystyle\frac{\la_{3,2}}{R_3}&0&\displaystyle\frac{\la_{3,4}}{R_3}&0&...&0&0\\
\vdots&\vdots&\vdots&\vdots&\vdots&\ddots&\vdots&\vdots\\
0&0&0&0&0&...&0&\displaystyle\frac{\la_{N-1,N}}{R_{N-1}}\\
0&0&0&0&0&...&\displaystyle\frac{\la_{N,N-1}}{R_N}&0
\end{pmatrix},
\]
where
\[
R_m=\begin{cases}\la_{m,m+1}+\la_{m,m-1}= 2\la m,\mbox{ if }m\in\{2,...,N-1\},\\
\la_{1,2}=\la,\mbox{ if }m=1,\\
\la_{N,N-1}=\la N,\mbox{ if }m=N.
\end{cases}
\]
The value of $R_m$ is the rate at which $X^N(t)$ leaves the state $m$. For the case $m\in\{2,...,N-1\}$    the transition probability is $\displaystyle\frac{\la_{m,m+1}}{R_m}=\displaystyle\frac{\la_{m,m-1}}{R_m}=\frac12$, and for $m=1$ and $m=N$ the transition probabilities are $\displaystyle{\frac\la{R_1}=\frac{\la N}{R_N}}=1$.
The two transition probabilities equal to 1 represent the reflection from  states 1 and $N$. 

Using the detailed balance equations, we obtain the stationary distribution $\wh\pi$ of  $\wh{X}(k)$
 \[
\wh\pi(1)=\wh\pi(N)=\frac1{2(N-1)},\ \wh\pi(m)=\frac1{N-1},\ m=2,...,N-1.
\] 
\begin{remark}
The chain $\wh{X}(k)$, with initial state $\wh{X}(0)$ is a simple random walk in $\mathcal{N}=\{1, \ldots, N \}$ with reflections on the boundaries. If the current state $\wh{X}(k-1) \in\{2,...,N-1\}$, then the next state $\wh{X}(k) = \wh{X}(k-1) + Y_k$, where the independent increments $Y_k$ take values $\pm 1$ uniformly $\PP(Y_k=1)= \PP(Y_k=-1)=1/2$, and $\PP(Y_k=1)=1\ \left( \PP(Y_k=-1)=1\right)$ on the boundary state, when $\wh{X}(k-1)=1\ \left( \wh{X}(k-1)=N \right)$.
\end{remark}

The dynamics of the process $(X^N(t))$ on continuous-time changes dramatically compared to the embedded chain $\wh{X}(k)$. The occupation (or sojourn) time plays a crucial role in the nature of the dynamics. Our aim is to  understand the occupation time behavior on different subsets of $\mathcal N$. 

For any $N$ the process $X^N(t)$ is homogeneous in time, but non-homogeneous in space. Thus, the sojourn time in any subset $\mathcal B\subset\mathcal N$ may vary depending on  $N$. We have shown that for any $M<N$ which is not very close to $N$  the time of the sojourn  probability in $\mathcal M=\{1, \ldots, M\}$  tends to 1 as $N, M\to \infty$ (see Theorem~\ref{22.03.23-1}). 

\section{Scaling} In this section we consider the scaled version of the processes $X^N$: $\xi^N(t)=\frac{X^N(t)}N$. The state space of $\xi^N$ is  $\A{N}^N=\left\{\frac1N,\frac2N,...,1\right\}$. The generator of $\xi^N(t)$ is the following operator (cf. \reff{3.1}) acting on the function space $\A{F}^N=\{{f}:\:\A{N}^N\to\R\}$ for all $\gamma\in\A{N}^N\subset [0,1]$
\bel{g4.1}
g^Nf(\mathscr{\gamma})=
\la N\gamma\Bigl(f\Bigl(\gamma+\frac1N\Bigr)-f(\gamma)\Bigr)\1_{\{\gamma<1\} }(\gamma)+\la N\gamma\Bigl(f\Bigl(\gamma-\frac1N\Bigr)-f(\gamma)\Bigr)\1_{\{\gamma > \frac1N\} }(\gamma).
\ee

Denote $\xi^N_{\gamma^N}$ the process $\xi^N$ with the initial point $\gamma^N\in\A{N}^N$, that is $\xi^N_{\gamma^N}(0)=\gamma^N$. Let $\mathbf D^{N}[0,T]$ be the set of all trajectories of $\xi^N$ during the time interval $[0,T]$ and let
\[
\mathbf D[0,T]=\bigcup_{N}\mathbf D^{N}[0,T].
\]
Since $\mathbf D^N[0,T]\subset\mathbf D[0,T]$ then all measures $\mu^N_{\gamma^N}$ of the prcesses $\xi^N_{\ga^N}$ are concentrated in $\mathbf D[0,T]$. However, the set $\mathbf D[0,T]$ does not contain all continuous functions. Let $\bar{\mathbf D}[0,T]$ be completion of $\mathbf D[0,T]$ with respect to the uniform metric. Then $\mathbf C[0,T]\subseteq\bar{\mathbf D}[0,T]$.

Let $(\gamma^N)$ be the sequence  of the initial points such that $\gamma^N\to \gamma\in(0,1]$, for example $\gamma^N=\frac{\lfloor \gamma N\rfloor}{N}$. 
In what follows we will show that the measures $\mu^N_{\gamma^N}$ converge to the measure $\mu_{\gamma}$ concentrated on the path set $\{\Ab{\gamma}:\:\Ab{\gamma}(0)=\gamma\}$. Further the bold letters we reserve for the functions. For example, $\Ab{\gamma}(t) = \gamma$ means that the function $\Ab{\gamma}$ at the point $t$ takes the value $\gamma$. We hope this notation will not confuse a reader.

\section{The Law of Large Numbers}\label{5.1}
  
In this section, we study the dynamics of a sequence of the processes $\xi^{N}$  as $N$ grows. 
One of the purposes of this study is to determine the evolutionary behavior of the scaled processes $\xi^{N}$ with the initial state from $\A{N}=\bigcup_N\A{N^{N}}$. 
Recall that the process $\xi^{N}_{\gamma^N}(t)$ denotes the process  $\xi^{N}(t)$ with the initial $\xi^N(0)=\gamma^N\in\A{N}$. The measure $\mu_{\gamma^N}$ which is the distribution of $\xi^{N}_{\gamma^N}(t)$ is concentrated in the path set $\mathbf D^{N}[0,T]$. 

We consider two cases for the law of large numbers: the single-point initial distribution (Section~\ref{4.1}) and the stationary initial distribution (Section~\ref{4.2}). 

 

\subsection{The single-point initial distribution}\label{4.1}

We introduce a set of constant functions. Let  $\Ab{\zeta}(t)\equiv \gamma \in [0,1]$, where $\gamma$ is a fixed number. Let $\gamma^N$ be a sequence of the initial states, $\xi^N_{\gamma^N}(0)=\gamma^N$, such that $\lim_{N\to\infty} \gamma^N = \gamma\in[0,1]$ (for example, we can choose $\gamma^N=\lfloor \gamma N\rfloor /N$). Here we prove the law of large numbers for $\xi^{N}_{\gamma^N}(t)$ when $N$ goes to infinity.  

\begin{theorem}
For any $\gamma\in [0,1]$ 
the relation
\bel{7.1}
\lim_{N\to\infty} \PP \left( \sup_{t\in[0,T]} \left|  \xi^{N}_{\gamma^N}(t) -  \Ab{\zeta}(t)\right| \ge \varepsilon \right) = 0.
\ee
holds true.
\end{theorem}
\proof Let $J^N(t)$ be the number of the jumps of $\xi^{N}_{\gamma^N}(t)$ on $(0,t]$. Consider an auxiliary process $Y(t)=\sum_{i=1}^{J^N(t)}Y_i,$ on the time interval $[0,T]$, where $Y_i$ are i.i.d. sequence of random variables with $P(Y_i=1)=1-P(Y_i=-1)=\frac 12$. Let $Y_N(t) := \frac 1N Y(t)$. 
Observe that the following inequality holds:
$$
\PP \left( \sup_{t\in[0,T]} \left| \gamma_N + Y_N(t) -\Ab{\zeta}(t)\right| < \varepsilon \right) \le
\PP \left( \sup_{t\in[0,T]} \left| \xi^{N}_{\gamma^N}(t) -  \Ab{\zeta}(t)\right| < \varepsilon \right).
$$
 Note that there exists $N_\varepsilon$ such that $|\gamma^N-\gamma|<\varepsilon$ for all $N>N_\varepsilon$. For these $N$ we have
$$
\PP \left( \sup_{t\in[0,T]} \left| Y_N(t) \right| < \frac{\varepsilon}{2} \right) <
\PP \left( \sup_{t\in[0,T]} \left| \gamma^N + Y_N(t) -\Ab{\zeta}(t)\right| < \varepsilon \right).
$$
To finish the proof of \reff{7.1} we use Kolmogorov inequality for the maximum of sum of i.i.d random variables (see, for example, \cite{Bil}, Theorem 22.4). We have
$$
\begin{aligned}
	\PP \left( \sup_{t\in[0,T]} \left| Y_N(t) \right| \ge \varepsilon \right)
	&= \sum_{n=0}^\infty
	\left. \PP \left( \sup_{t\in[0,T]} \Bigl| \sum_{i=1}^{J^N(t)} Y_i \Bigr| \ge \varepsilon N \ \right| \  J^N(T) = n \right)  \PP \left( J^N(T)=n \right) \\
	& = \sum_{n=0}^\infty
	\left. \PP \left( \max_{1 \le k \le n } \Bigl| \sum_{i=1}^{k} Y_i \Bigr| \ge \varepsilon N \ \right| \  J^N(T) = n \right)  \PP\left( J^N(T)=n \right) \\
	& \le \sum_{n=0}^\infty  \frac{n}{4\varepsilon^2 N^2} \PP(J^N(T)=n) =
	{ \frac{1}{4\varepsilon^2 N^2} \mathbb{E}(J^N(T)) 
		\leq \frac{ N T }{4\varepsilon^2 N^2} }
	\to 0 \mbox{ as } N\to \infty.
\end{aligned}
$$
Finally,
\begin{equation}\label{eq2}
\begin{aligned}
	1 & \ge \PP \left( \sup_{t\in[0,T]} \left| \xi^{N}_{\gamma^N}(t) -  \Ab{\zeta}(t)\right| < \varepsilon \right) \ge 1- \PP \left( \sup_{t\in[0,T]} \left| Y_N(t) \right| \ge \frac{\varepsilon}{2} \right) \\
	& \ge 1- \frac{ T }{\varepsilon^2 N} \to 1,  \mbox{ as } N\to \infty.
\end{aligned}
\end{equation}
This finishes the proof of \reff{7.1}. $\Box$
\begin{remark}
We apply the inequality $\mathbb{E}(J^N(T))\leq NT $ because the intensity of the jumps is maximal near the boundary on 1.
\end{remark}

\subsection{Stationary initial distribution}\label{4.2}
We now turn to the case when the initial state is distributed according to the stationary measure. To this end, we introduce the constant process $\Ab{\zeta}(t)$ as before. That is $\Ab{\zeta}(t)\equiv \gamma\in[0,1]$ for $t\in[0,T]$ if $\Ab{\zeta}(0)=\gamma$. The distribution of $\Ab{\zeta}(t)$ satisfies the following conditions: for any $u\in(0,1]$ the probabilities $\PP(\Ab{\zeta}(t)<u)=1$ and $\PP(\Ab{\zeta}(t)>u)=0$. The described property means that the process $\Ab{\zeta}(t)$ is concentrated on the line $\Ab{\ga}(t)\equiv 0$.

Let $\xi_{\mathbf{st}}^N(t)$ be the piecewise constant process on $\mathbf{D}^{N}[0,T]$ with the initial distribution $\pi^N$. Note that the distribution $\pi^N$ is defined on the set $\mathcal N$, but we will maintain the same notation for the stationary measure on $\A{N}^N$. Then the distribution $\xi_{\mathbf{st}}^N(t)$ at any $t\in[0,T]$ is also $\pi^N$ since $\pi^N$ is  stationary. Recall that the generator of $\xi_{\mathbf{st}}^N(t)$ is \reff{g4.1}.

\begin{theorem}\label{t4.4}
\bel{7.0}
\xi^{N}_{\mathbf{st}}\Rightarrow\Ab{\zeta}\mbox{ \ weakly}\ee
as $N\to\infty$
\end{theorem}

 The proof of the theorem is based on the result of the paper \cite{BP} on the convergence of integrals of finite-dimensional distributions of random processes. The paper \cite{BP} states conditions for the weak convergence of the functionals
 $
 \Phi(\xi_{st}^N)=\int_0^T\phi(\xi_{\bold{st}}^N(u)\ed u
$ to  $\Phi(\Ab{\zeta})=\int_0^T\phi(\Ab{\zeta}(u))\ed u$, where  $\phi$ is the continuous function from $\mathbf C[0,1]\to \R$.

The mentioned theorem in the considered case is 
\begin{theorem}[Borovkov, Pechersky, \cite{BP}]\label{BPtheorem}
	Let 
	\[
	\lim_{N\to\infty}\int_0^T\dots\int_0^T\PP\left(\bigcap_{i=1}^\ell(\xi^{N}_{\bold{st}}(t_i)<u_i)\right)\ed t_1,...,\ed t_\ell=\int_0^T\dots\int_0^T\PP\left(\bigcap_{i=1}^\ell(\Ab{\zeta}(t_i)<u_i)\right)\ed t_1,...,\ed t_\ell
	\]
	for any finite number $t_1,...,t_\ell\in[0,T]$ then the distribution of $\Phi(\xi^{N}_{\bold{st}})=\int_0^T\phi(\xi^{N}_{\bold{st}}(u))\ed u$ weakly converges to $\Phi(\Ab{\zeta})=\int_0^T\phi(\Ab{\zeta}(u))\ed u$ for any $\phi\in\mathbf C[0,1]$.
\end{theorem}

Due to the properties of the process $\Ab{\zeta}$, the convergence of the integral functionals is sufficient for weak convergence \reff{7.0}. The result we  apply to prove \reff{7.0} reformulated for our case as the following  
\begin{lemma}\label{lemma} 
For any $ \{u_1,...,u_\ell\}\subset [0,1]$ and $\Ab{\zeta}<\wt{u}=\min_{j=1,...,\ell}u_j$ the following equality
\bel{7.2}\lim_{N\to\infty}\int_0^T\dots\int_0^T\PP\left(\bigcap_{i=1}^\ell(\xi^{N}_{\bold{st}}(t_i)<u_i)\right)\ed t_1,...,\ed t_\ell=\int_0^T\dots\int_0^T\PP\left(\bigcap_{i=1}^\ell(\Ab{\zeta}(t_i)<u_i)\right)\ed t_1,...,\ed t_\ell=T^\ell
\ee
holds for any finite set of time points $t_1,...,t_\ell\in[0,T]$.
\end{lemma}

\noindent
\Proof \textit{of Lemma \ref{lemma}}. The righthand side of \reff{7.2}   
$$
\int_0^T  \dots\int_0^T  \PP\left( \cap_{j=1}^\ell (\Ab{\zeta}(t_j)< u_j)\right) dt_1 \ldots dt_\ell =T^\ell,
$$
by the definition of $\Ab{\zeta}(t)$.
To find the limit on the lefthand side, let $0<a < \wt{u} = \min_{j=1,...,\ell}\{u_j\}$ and  $\e>0$ be such that $(a-\e, a+\e) \subset(0,\wt{u})$. For a fixed $N$ we have
\begin{figure}[h]
	\caption{Illustration for the proof of Theorem~\ref{t4.4}}	
	\centering 
	\includegraphics[width=15cm]{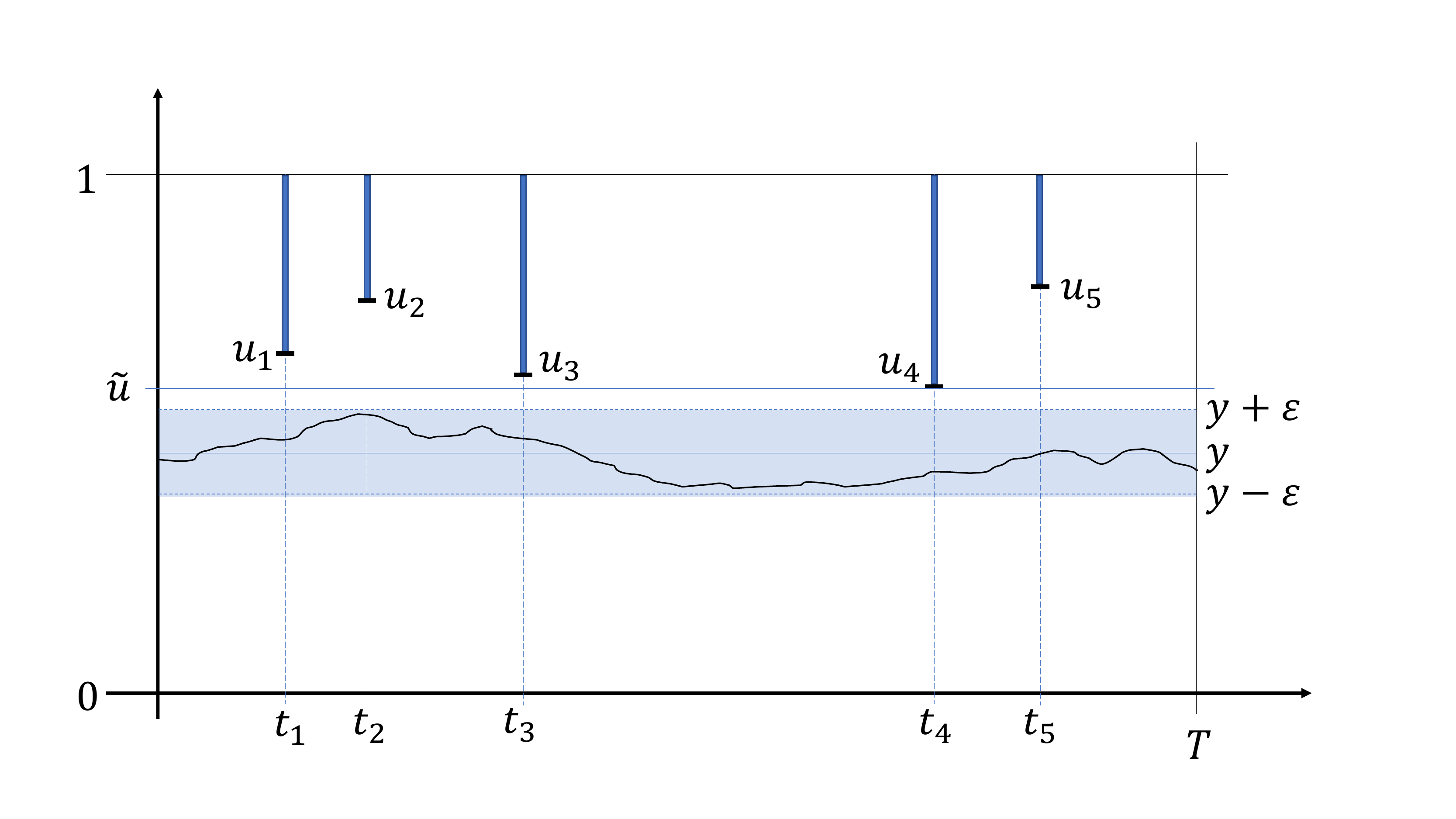} 
	\label{fig1}
\end{figure}

\begin{equation}\label{eq1}
	\begin{aligned}
		\int_0^T  \dots\int_0^T & \PP\left( \cap_{j=1}^\ell \big( \xi^N_{\mathbf{st}}(t_j)< u_j\big)\right)dt_1 \ldots dt_\ell 
		\\
		&=
		\int_0^T  \dots\int_0^T\sum_{ y < a }\PP\left( \cap_{j=1}^\ell \big( \xi^N_y(t_j)< u_j\big)\right)\pi^N(y)dt_1 \ldots dt_\ell 
		\\
		& + 
		\int_0^T  \dots\int_0^T\sum_{ y \ge a }\PP\left( \cap_{j=1}^\ell \big( \xi^N_y(t_j)< u_j\big)\right)\pi^N(y)dt_1 \ldots dt_\ell 
	\end{aligned}
\end{equation}
Probabilities in the last two sums we estimate as the following. For any $y < a$, and sufficiently large $N$, we use the following bounds (see \reff{eq2})
$$
1 \ge \PP\left( \cap_{j=1}^\ell \big( \xi^N_y(t_j)< u_j\big)\right) \ge \PP \left(\xi^N_y(t)\in [y-\e,y+\e]\right) \ge 1 - \frac{T}{\varepsilon^2 N}.
$$
For the first term of \reff{eq1} we have
$$
T^\ell \ge \int_0^T  \dots\int_0^T\sum_{ y < a }\PP\left( \cap_{j=1}^\ell \big( \xi^N_y(t_j)< u_j\big)\right)\pi^N(y)dt_1 \ldots dt_\ell  \ge T^\ell \pi^N([0,a)) + O\left(\frac1N\right) \to T^\ell,
$$
as $N\to\infty$, according to Theorem~\ref{4.0}.  In the second term for all $y \ge x$ corresponding probabilities we can bound by 1, and
$$
0 \le \int_0^T  \dots\int_0^T\sum_{ y \ge a }\PP\left( \cap_{j=1}^\ell \big( \xi^N_y(t_j)< u_j\big)\right)\pi^N(y)dt_1 \ldots dt_\ell
\le 
T^\ell \pi^N([a,1]) \to 0,
$$
as $N\to\infty$, according to Theorem~\ref{4.0}.
This finishes the proof of the lemma. $\Box$

\noindent
\Proof \textit{of Theorem~\ref{t4.4}}. Theorem~\ref{BPtheorem} and Lemma~\ref{lemma} prove the theorem. $\Box$

 \section{Large Deviations} This section is devoted to studying the process deviations far from average values. This study belongs to the well-developed area of probability theory known as the theory of \textit{Large Deviations}. Here we use the techniques developed in \cite{Kurtz2}. We apply and adapt these results with additional calculations corresponding to the peculiarities of our model. We also refer a reader to \cite{K1} and \cite{Ber}, where the method was applied to some similar Markov processes. Thus, here we only pass through some points related to our process, and we would like to avoid cumbersome constructions used in the proof of LDP.
 
 \begin{remark}
 	The one-dimensional case of the generalized Ehrenfest model studied in \cite{K1} is very similar to the model studied here with a slight difference. Despite the differences in the generators of Fleming semigroups the corresponding Hamiltonians are similar. We will refer to this paper passing through the proof of LDP.	
 \end{remark}

Let $\mathbf C_1[0,T]$ be a set of continuously differentiable function on $[0,T]$ with values from $[0,1]$. 
The following two terms are our main interest.
\begin{enumerate}
\item[] \textit{Rate function (functional)} $I: \Ab{\ga}\to \mathbb{R}$ for $\Ab{\ga}\in \mathbf C_1[0,T]$, see \reff{RF}. 

\item[] \textit{Hamiltonian} $H(\gamma,\k)$,
with $\gamma\in [0,1]$ and $\k\in \mathbb{R}$, see \reff{11.2}. 
\end{enumerate}

For any $\Ab{\ga}\in \mathbf C_1[0,T]$ the rate function $I(\Ab{\ga})$ measures deviations of $\Ab{\ga}$ from the process average. $I(\Ab{\gamma})$ is approximately the logarithm of the following probability of the deviation $\Ab{\gamma}(t)-\mathbf E\xi^N(t)$:
 \[
 I(\Ab{\gamma})=\lim_{\varepsilon}\lim_N\frac1N\ln\PP\left( \xi^N(\cdot)\in U_\e(\Ab{\gamma}) \right),
 \] 
 where $U_\e(\Ab{\gamma})$ is $\e$-neighborhood of $\Ab{\gamma}$.
 If $\Ab{\gamma}(t)=\lim_N\mathbf E\xi^N(t)$ then $I(\Ab{\gamma})=0$. 
 The existence of this limit is known as a \textit{local large deviation principle (LLDP)} 
 
 \begin{definition} \label{d1.1}
 	Let the sequence of random processes $\xi^N(\,\cdot\,)$ be
 	defined as in \eqref{4.1}.
 	We say that this sequence satisfies a {\rm local large deviation principle}
 	with an {\rm{LD}} functional
 	$I\,:\mathbf{C}_1[0,T]\to [0,\infty ]$ if for all $f \in \mathbf{C}_1[0,T]$
 	$$
 	\lim\limits_{\varepsilon\rightarrow 0}\limsup\limits_{N\rightarrow \infty}\dfrac{1}{N}
 	\ln\mathbf{P}(\xi^N(\,\cdot\,)\in U_\varepsilon(f))
 	=\lim\limits_{\varepsilon\rightarrow 0}\liminf\limits_{N\rightarrow \infty}\dfrac{1}{N}
 	\ln\mathbf{P}(\xi^N(\,\cdot\,)\in U_\varepsilon(f))=-I(f),
 	$$
 	where
 	$$U_\varepsilon(f)=\{g\in \mathbf{D}[0,T]: \ \rho(f,g)<\varepsilon\}, \ \mbox{ and }
 	\ \rho(f,g) := \sup_{t\in [0,T]} \left| f(t) - g(t) \right|.$$
 \end{definition}
 
Frequently, the existence of this limit together with exponential tightness verify the \textit{large deviation principle (LDP)}.
 
\begin{definition} \label{d1.2} Let $\xi^N(\,\cdot\,)$ be sequence of random processes defined as in \eqref{4.1}. We say that it satisfies a large deviation principle with an {\rm{LD}} functional $I\,:\;\mathbf{C}_1[0,T]\to (0,\infty ]$ if, whenever measurable set $\mathbb{B}
	\subseteq\mathbf{C}_1[0,T]$, we have that
	$$\begin{array}{c} \limsup\limits_{T \rightarrow \infty} \dfrac{1}{N} \ln \mathbf{P}(\, \xi^N(\,\cdot\,)
		\in\mathbb{B} \,) \leq - I({\rm{Cl}}(\mathbb{B})),\\
		\liminf\limits_{T \rightarrow \infty} \dfrac{1}{N} \ln \mathbf{P}(\, \xi^N(\,\cdot\,)
		\in\mathbb{B}\,)\geq -I({\rm{Int}}(\mathbb{B})).\end{array}$$
\end{definition}
 
 This way worked well in papers \cite{VLSY3} -- \cite{MPY}, where the method to obtain the LLDP was proposed and applied for birth-and-death processes. But we couldn't extend the method for the process considered here, so we opted for the operator method proposed by Fend \& Kurtz \cite{Kurtz2}.

  It seems the shortest way to find the rate function is Bryc's theorem. It claims the following.
  \bel{10.1}
  I(\Ab{\gamma})=\sup_{\mathscr{L}}\left\{\mathscr{L}(\Ab{\gamma})-\Lambda(\mathscr{L})\right\},
  \ee
  where $\{\mathscr L\}$ is a set of  functionals on the process path set $\mathbf D[0,T]$, and 
  $$
  \Lambda(\mathscr{L})=\lim_N\frac1N\ln\int\exp\{\mathscr L(\Ab{\gamma})\} \PP(\xi^N\in\ed \Ab{\gamma})
  $$ 
 for any  functional $\mathscr{L}$. We have omitted the conditions and other details required for  Bryc's theorem  (see Proposition 3.8 in \cite{Kurtz2}).
  
Unfortunately, this is not a very practical way. The main difficulty is related to the fact that we consider Markov-dependent random variables. In this case, $\Lambda(\mathscr{L})$ is too general, and therefore, it is very rarely computable. Do the functionals $\Lambda(\mathscr{L})$ form a semigroup acting on a set of functions? A positive answer to this question would construct a generator of this semigroup, which makes it possible to create a research method of $\Lambda(\mathscr{L})$.

\subsection{Fleming semigroup}\label{S5.1} The above question  is solved with the help of a (Fleming) semigroup, proposed by Fleming (\cite{Fl}):
\bel{10.1}
V_t^N{f}(\gamma)=\frac1N\ln\int\exp\{N{f}(\ga)\}  \PP(\xi^N(t)\in \ed {\ga}) \mid \xi^N(0)=\gamma_0),
\ee
where the probe functions belongs to the set of $C_1[0,1]$ continuously differentiable function on $[0,1]$ with values from $\mathbb R$, $f\in\mathbf C_1[0,1]$.
The problem now is to investigate the limiting semigroup
\bel{10.2}
\mathbf V_t{f}(\gamma)=\lim_NV_t^N{f}(\gamma).
\ee
Since $V_t^N$ is a semigroup it is reasonable to find its generator. The generator corresponding to  $V_t^N$ describes the  Markov dynamics which is called  \textit{non-linear} and its generator is  called \textit{non-linear Hamiltonian $\mathcal H^N$}: if $\gamma \in\A{N}^N\setminus\{\frac1N, 1\}$, then
\[
 \begin{aligned}
 \mathcal H^N{f}(\gamma) &=\frac1N\exp\big\{-N{f}(\gamma)\big\}\left(\la N\gamma\left(\exp\left\{N{f}\left(\gamma+\tfrac1N\right)\right\}-\exp\big\{N{f}(\gamma)\big\}\right)\right. \\
 &\phantom{= \frac1N\exp\big\{-Nf(x_N)\big\}}
  + \la N\gamma\left(\exp\left\{N{f}\left(\gamma-\tfrac1N\right)\right\}-\exp\big\{N{f}(\gamma)\big\}\right)\Big),
  \end{aligned}
\]
if $\gamma = \frac1N$, then
$
\mathcal H^N{f}(\tfrac1N) = \frac{\lambda}{N}\exp\left\{-N{f}(\tfrac1N)\right\} \left(\exp\left\{N{f}\left(\tfrac2N\right)\right\}-\exp\left\{N{f}(\tfrac1N)\right\}\right),
$
and if $\gamma = 1$, then
$
\mathcal H^N{f}(1) = \lambda \exp\left\{-N{f}(1)\right\} \left(\exp\left\{N{f}\left(1-\tfrac1N\right)\right\}-\exp\left\{N{f}(1)\right\}\right).
$ 
It is easy to see that the limiting operator $\mathcal H$ should be
\be\label{11.3}
\mathcal Hf(\gamma) =\lim_N\mathcal H^N{f}(\gamma) =\la \gamma\left(\exp\left\{\frac{\ed f}{\ed \gamma}(\gamma)\right\}-1\right)+\la \gamma\left(\exp\left\{-\frac{\ed \gamma}{\ed \gamma}(\gamma)\right\}-1\right) =: H\left(\gamma, \frac{\ed f}{\ed \gamma} (\gamma) \right),
\ee 
where the set of probe functions is the set $\mathbf C_1^0[0,1]$ of continuously differentiable functions with zero derivative at 1, $\mathbf C_1^0[0,1] \subset \mathbf C_1[0,1]$. Using the notation $\k=\frac{\ed f} {\ed \gamma}(\gamma)$ we write the Hamiltonian 
\bel{11.2}
H(\gamma,\k)=\la \gamma(e^\k-1)+\la \gamma(e^{-\k}-1).
\ee

The proof of the operator convergence 
\[
\mathcal H=\lim\mathcal H^N
\]
in the $\mathbf C_1^0[0,1]$ is not diffcult. First, we consider the convergence on $\A{N}$, where the limit 
\bel{11.1}
\lim_{N\to\infty} \sup_{\gamma\in \A{N}} \left| \mathcal H^N{f}(\gamma) - \mathcal Hf(\gamma) \right| =0,
\ee
is obvious.
Second, consider a sequence $(\gamma_N)$ where $\gamma_N\in\A{N}$ such that  $\gamma_N \to \gamma$ for an arbitrary  $\gamma\in \mathbf [0,1]$. Then 
$$
\lim_N\frac{f(\gamma_N + 1/N) - f(\gamma_N)}{1/N} = \frac{\ed f}{\ed \gamma} (\gamma)
$$
and the derivative of $f$ is continuous on $[0,1]$. Therefore \reff{11.1} can be extended to
\[
\lim_{N\to\infty} \sup_{\gamma\in [0,1]} \left| \mathcal H^N{f}(\gamma) - \mathcal Hf(\gamma) \right| =0.
\]

As can be seen the resulting function $\mathcal H f(\gamma)=H\left(\gamma, \frac{\ed {f}} {\ed \gamma}\right)$ depends on the derivative $\frac{\ed {f}}{\ed \gamma}$, the first variable $\gamma$ is the reason for the non-linearity of the dynamics.

Next considerations prove that $\mathcal H$ is the generator of the semigroup $\mathbf V_t$.

\subsection{Exponential compact containment}\label{S5.2}
The {\it exponential compact containment } can be obtained directly from the definition for our case. The sequence $\xi^N$ satisfy exponential compact containment condition in $\bar{\mathbf D}(0,T]$ iff for each $a>0$ there exists a compact $K_a\subset [0,1]$ such that 
$$
\limsup_{n\to\infty} \frac{1}{n} \ln P( X_n(t) \notin K_a \mbox{ for some }t\le T) \le - a.
$$ 
In our case this condition is easily verified taking $K_a = [0,1]$. Observe that it is not exponential tightness. But together with generator convergence and the comparison principle (see further for the definition) exponential containment will imply exponential tightness according to Theorem 6.14, \cite{Kurtz2}.

\subsection{Fleming semigroup generator}\label{S5.3} We now  consider when $\mathcal H$, \reff{11.3}, is the generator of the Fleming semigroup $\mathbf V_t$, \reff{10.2}. The proof is a rather  cumbersome sequence of arguments. One of the important conditions for $\mathcal H$ to be the generator is the existence of solutions to the equation  $(I-\beta \mathcal H)f=h$ for each $h$ and small positive $\beta$ (Chapter 6, \cite{Kurtz2}). 
It is so called the \textit{viscosity equation}. This condition is difficult to verify. The \textit{comparison principle} overcomes this difficulty. 

Let us rewrite the equation using Hamiltoniam:
\be\label{HamEq}
f(\gamma) - \beta H\left( \gamma, \tfrac{\ed f}{\ed \gamma} (\gamma)\right) - h(\gamma) =0.
\ee
Let $u$ and $v$ be subsolution and supersolution of equation \reff{HamEq} (see Definition~\ref{subsupersolutions} in Appendix 1). We say that the equation verifies the 
\textit{comparison principle} if for any subsolution $u$ and supersolution $v$, it is verified that $u\le v$.

The proof of the comparison principle for the case under study is given in Appendix 1, see Proposition~\ref{CompP}. When the comparison principle holds, if there exists a viscosity solution, it is unique.

\subsection{Rate function. Variational representation.}\label{S5.4}

We have established that Hamiltonian $H$ is Fleming group generator. The next our goal is to represent the rate function $I(\Ab{\ga})$ as Legendre transform of the Hamiltonian that is to  obtain the \textit{variational form} of the rate function
\bel{16.1}
I(\Ab{\gamma}) =   
\int_0^T \sup_{\k} \left\{ \k\dot{\Ab{\gamma}}(t) - H(\Ab{\gamma}(t),\k) \right\}\ed t.
\ee
This expression is difficult to obtain directly from Fleming semigroup terms (see \reff{10.1}, \reff{10.2}) and require some results from control theory. The Fleming semigroup in the large deviation is interpreted as Nisio semigroup in control theory.
The proof of \reff{16.1} uses a method developed in a number of works.  
The detail presentation of the method is represented in Capter 8 from \cite{Kurtz2}, see also \cite{K1}, \cite{Ber}, for the application of the method. 
 
 The method is a series of transformations of the generator $g^N$ \reff{4.1}. The first non-linear transformation provided the Fleming semigroup $\mathbf{V}_t$ and its generator is $H= \lambda \gamma \left( e^{\k} + e^{-\k} -2 \right)$, see Section~\ref{S5.1} and Section~\ref{S5.3}. 
 The Fleming semigroup has the exponential form not easy to investigate. The last transformation is the Nisio semigroup construction. 

\subsubsection{Nisio Semigroup} In this stage the method shows that the Fleming semigroup coincides with the so-called \textit{Nisio semigroup} $\mathbf{W}(t)$ from the control theory (see \cite{Kurtz2}, Chapter 8). The control  theory provides the variational representation \reff{16.1} for the rate function.

Before starting let us observe that $H(\gamma,\k) $ is a convex function on the second coordinate $\k$, and 
\begin{equation}\label{Hxp}
H(\gamma,\k) = \sup_{\gamma} \left( \k \gamma - L(\gamma,\k) \right),
\end{equation}
where 
\begin{equation}\label{Lxu}
L(\gamma,\k) :=\sup_{\k} \left( \k\gamma - H(\gamma,\k) \right).
\end{equation}
The Nisio semigroup has the following form
\begin{equation}\label{Nisio1}
	\mathbf{W}(t)f(x) = \sup_{\Ab{\gamma} \in \mathcal{AC}[0,T], \\ \Ab{\gamma}(0)=\gamma} f(\gamma(t)) - \int_0^t L(\gamma(s),\dot{\gamma}(s)) ds,
\end{equation}
where $\Ab{\ga}\in\mathcal{AC}[0,T]$, and $\mathcal{AC}[0,T]$ is the set of all absolutely continuous paths on $[0,T]$. 

A more easy approach is to use a \textit{relaxed presentation} of Nisio semigroup instead of \reff{Nisio1}.  Let $U=\{\dot{\Ab{\gamma}}(\cdot)\}$ and $\mathcal{M}_m(U)$ be the space of Borel measures $\lambdabar$ on $U\times[0, T]$ satisfying $\lambdabar (U \times [0, t]) = t$ for all $t\in [0, T]$. The measure $\lambdabar$ is known as a {\rm relaxed control} (see detailed definition in Appendix 2). 

The relaxed controls determine the Nisio semigroup as
\begin{equation}\label{Nisio2}
	\mathbf{\bar{W}}(t)f(x) = \sup_{(\Ab{\gamma},\lambdabar),\ \Ab{\gamma}(0)=x} f(\Ab{\gamma}(t)) - \iint_{U\times[0,1]} L(\Ab{\gamma}(s),u) \lambdabar(du \times ds),
\end{equation}

\subsubsection{Variational representation of rate function}
Finally, the relation \reff{16.1} follows from \cite{Ber}, Theorem 3.1, which states that
	\begin{enumerate}
		\item $\mathbf{V}_t(f) = \bar{\mathbf{W}}(t)f$ 
		\item $I({\Ab{\gamma}}) = \inf_{\lambdabar:\ (\Ab{\gamma},\lambdabar)} \iint_{U\times [0,T]} L(\Ab{\gamma}(s),u)\lambdabar(du\times ds)$;
		\item the rate function can be written as an action integral:
		\begin{equation}\label{21.03.23-1}
			I(\Ab{\gamma}) = \int_{0}^T L(\Ab{\gamma}(s),\dot{\Ab{\gamma}}(s))ds,
		\end{equation}
	if $\Ab{\gamma}$ is absolutely continuous function. 
	\end{enumerate}
We have omitted the required conditions. 

\subsection{Hamiltonian system} 
At that point, we have made four steps to obtain LDP and the rate function. 

\begin{enumerate}
\item[] 
Step 1. Section~\ref{S5.1} established the convergence of the sequence of operators $\mathcal{H}^N$, where we obtained the limiting operator $\mathcal{H}$.

\item[]
Step 2. Section~\ref{S5.2} checked the exponential compact containment condition. 

\item[]
Step 3. Section~\ref{S5.3} with Appendix 1 established the comparison principle. It proves that the limiting operator $\mathcal{H}$ generates the semigroup $\mathbf{V}_t$: $\frac{d}{dt} \mathbf{V}(f)=\mathcal{H}(\mathbf{V}_t(f))$.
\end{enumerate} 

\noindent
According to Theorem 6.14 from \cite{Kurtz2} these three steps assure that the sequence of our processes is exponentially tight and satisfies an LDP with rate function expressed via the Fleming semigroup (see formulas for rate functions (5.30) and (5.31) in \cite{Kurtz2}). Practically the exponetial tightness implies that the rate function is a good rate function (see, for example, Theorem 3.7 \cite{Kurtz2}).

This result does not provide a useful presentation for the rate function. The next fourth step of Feng \& Kurtz theory provides the practical version of the rate function based on control theory.

\begin{enumerate}
	\item[] Step 4. Section~\ref{S5.4} with Appendix 2 state and checked conditions when the rate function $I$ can be written as an action integral \reff{21.03.23-1}. (Condition 8.9 and 8.11 from \cite{Kurtz2} are stated in Appendix 2. They are easily verified. See, for example, Propositions 3.6 and 3.8 from \cite{Ber}). 
\end{enumerate}

Here, based on a variational representation of the rate function, fortunately, we can obtain the rate functional \reff{RF} explicitly. Moreover, we can find the optimal path for the following rare event: starting at the point $\Ab{\gamma}(0)$ the process will stay at time $T$ at the point $\Ab{\gamma}(T)$ far from their mean, $\Ab{\gamma}(0) \ne \Ab{\gamma}(T)$ (see example in Appendix 3.)

Let us now find paths between fixed boundary conditions. To this end, we should investigate the following Hamiltonian system
\[
\left\{
\begin{aligned}
	&\dot {\Ab{\gamma}}=\la \Ab{\gamma}
	\left(e^\k-e^{-\k}\right)\\
	&\dot \k=-\la \left(e^\k+e^{-\k}-2\right)
\end{aligned}
\right.
\]
First equation of Hamiltonian system is the derivative on $\k$,  $\dot{{\Ab{\gamma}}} = \frac{\partial H}{\partial\k}$, to find the supremum. The second one is Euler–Lagrange equation.

\vspace{0.5cm}
\noindent
\textbf{Rate function}: The first equation provides the explicit formula for the rate function. Let $z=e^{\kappa}$. For given functions $\Ab{\gamma}(t), \dot{\Ab{\gamma}}(t)$ we obtain the quadratic equation on $z$
$$
z^2 - \frac{\dot{\Ab{\gamma}}}{\lambda \Ab{\gamma}} z -1 =0,
$$
with admissible solution 
$$
z(t)=e^{\kappa(t)} = \frac{\dot{\Ab{\gamma}}}{2 \lambda \Ab{\gamma}} + \sqrt{\left( \frac{\dot{\Ab{\gamma}}}{2\lambda \Ab{\gamma}} \right)^2 +1}.
$$
Finally,
\be\label{RF}
I(\Ab{\gamma}) = \int_0^T  \left[ \dot{\Ab{\gamma}} \ln \left( \frac{\dot{\Ab{\gamma}}}{2 \lambda \Ab{\gamma}} + \sqrt{\left( \frac{\dot{\Ab{\gamma}}}{2\lambda \Ab{\gamma}} \right)^2 +1} \right) - \dot{\Ab{\gamma}} + 2\lambda \Ab{\gamma} - \frac{(2\la\Ab{\gamma})^2}{\dot{\Ab{\gamma}} + \sqrt{\dot{\Ab{\gamma}}^2 + (2\la \Ab{\gamma})^2}} \right] dt
\ee

\vspace{0.5cm}
\noindent
\textbf{Solution of the system}:
Let $z=e^\k$ then $\dot\k=\frac{\dot z}z$. The second equation of the above system is 
\[
\frac{\dot z}z=-\la\left( z+\frac1z\right)+2\la,
\]
or
\[
\dot z=-\la(z-1)^2.
\]
The solution is
\[
z(t)= \left\{ 
\begin{aligned}
	\frac1{\la t - c_1}+1, & \mbox{ when }\Ab{\gamma}(t) \nequiv const; \\
	1, & \mbox{ when }\Ab{\gamma}(t) \equiv const.
\end{aligned} 
\right.
\]
where $c_1$ is some constant.   

Now the first equation of the system (when $\Ab{\gamma} \nequiv const$) is
\[
\dot{\Ab{\gamma}}= \lambda \Ab{\gamma} \left( z(t) - \frac1{z(t)} \right) = \la\Ab{\gamma}\left(\frac1{\la t - c_1} + \frac1{\la t - c_1+1}\right),
\]
or
\[
\frac{\ed\  \ln(\Ab{\gamma})}{\ed t} =\la \left(\frac1{\la t - c_1} + \frac1{\la t - c_1+1}\right)
\]
Then
\[\Ab{\gamma}(t)=
\left\{ 
\begin{array}{rl}
	(\la t - c_1)(\la t - c_1 +1)c_2, & \mbox{ if }\Ab{\gamma}(t)\nequiv const \\
	const, & \mbox{ if }\Ab{\gamma}(t)\equiv const \in [0,1].
\end{array}
\right.
\]
where $c_2$ is some positive constant. The parameters $c_1, c_2$ are uniquely determined by the boundary conditions $\Ab{\gamma}(0), \Ab{\gamma}(T)$ and by the fact that the solution $\Ab{\gamma}(t)\in [0,1]$ when $t\in [0,T]$. 

The examples of the solutions are in Appendix 3.

\section*{Acknowledgment}
AY thanks A. Logachov for fruitful discussions. AY thanks FAPESP for the financial support via the grant 2017/10555-0.

\section*{Appendix 1: comparison principle}\label{A33}

The following proposition proves the comparison principle. It is similar to Proposition 3 in \cite{K1}. We adapt the proof in \cite{K1} for our case. Let us recall some basic definitions from \cite{K1} slightly adapted for our notations. Recall, see \reff{11.2}, that we are working with our Hamiltonian $H(x,\k) = \lambda x \left( e^{\k} + e^{-\k} -2 \right)$.

\begin{definition}\label{subsupersolutions}
	$u$ is called a (viscosity) {\rm subsolution} if
	\begin{enumerate}
		\item[] $u$ is bounded, upper semi-continuous; 
		\item[] for every $\phi \in \mathbf{C}_1[0,1]$ and $x_0:=\mbox{argmax}(u(x) - \phi(x))$ we have $u(x_0) - \beta H(x_0,\phi^\prime (x_0)) - h(x_0) \le 0.$ 
	\end{enumerate}
	We say that $v$ is a (viscosity) {\rm supersolution} if
	\begin{enumerate}
		\item[] $v$ is bounded, lower semi-continuous; 
		\item[] for every $\phi \in \mathbf{C}_1[0,1]$ and $x_0:=\mbox{argmin}(v(x) - \phi(x))$ we have $v(x_0) - \beta H(x_0,\phi^\prime (x_0)) - h(x_0) \ge 0.$ 
	\end{enumerate}
	We say that $u$ is a (viscosity) \textit{solution} if it is both a sub and super solution.
\end{definition}

\begin{proposition}\label{CompP}
	 Let $\beta>0$ and $h\in \mathbf{C}[0,1]$. Then the comparison principle holds true for 
	\begin{equation}\label{comp1}
		f(x) - \beta H(x,f^\prime (x)) - h(x) = 0.
	\end{equation}
\end{proposition}

\noindent
{\it Proof.} With small adaptation the proof follows \cite{K1}, Section 4. See especially \cite{K1}, Section 4.1, where comparison principle verified for a one-dimensional Erenfest model similar to our case. 

For given $\beta>0$, $h\in \mathbf{C}[0,1]$ let $u, v$ be some sub- and super-solutions of \reff{comp1}. Choose the following penalized function $\Psi(x,y) = \frac{1}{2} (x-y)^2$. For a positive increasing parameter $\alpha$ (for example, we can suppose $\alpha$ as natural numbers) we define the following points $x_\alpha, y_\alpha \in [0,1]$
$$
u(x_\alpha) - v(y_\alpha) - \frac{\alpha}{2}(x_\alpha-y_\alpha)^2 = \sup_{x,y\in E} \left\{ u(x) - v(y) - \frac{\alpha}{2}(x-y)^2\right\}. 
$$
Since the functions $u, v$ are bounded the limit $\lim_{\alpha\to\infty} \alpha \Psi(x_\alpha, y_\alpha)$ should be zero. All limit points of $(x_\alpha, y_\alpha)$ should take the form $(z,z)$  such that 
$$
u(z) - v(z) = \sup_{x\in E} u(x) - v(x).
$$
It is the general result, and we refer to Lemma 4 in \cite{K1}. 

For such chosen sequence of points $x_\alpha, y_\alpha$ according to Proposition 2 in \cite{K1} we need to prove that
$$
\liminf_{\alpha \to \infty} H(x_\alpha, \alpha(x_\alpha-y_\alpha)) - H(y_\alpha, \alpha(x_\alpha-y_\alpha)) \le 0.
$$
In our case
\be\label{CPe1}
H(x_\alpha, \alpha(x_\alpha-y_\alpha)) - H(y_\alpha, \alpha(x_\alpha-y_\alpha)) = \lambda(x_\alpha - y_\alpha) \left( e^{\alpha(x_\alpha-y_\alpha)} + e^{-\alpha(x_\alpha-y_\alpha)} -2\right).
\ee
Lemma 5 in \cite{K1} provides more control: accoding this lemma we have
\be\label{CPe2}
\sup_{\alpha} H(y_\alpha, \alpha (x_\alpha - y_\alpha)) = 
\lambda y_\alpha \left( e^{\alpha(x_\alpha-y_\alpha)} + e^{-\alpha(x_\alpha-y_\alpha)} -2\right) < \infty.
\ee
Thus, suppose that $z\in (0,1]$ then, according to \reff{CPe2} the sequence $\alpha(x_\alpha - y_\alpha)$ is bounded, and it provides \ref{CPe1} (see some details in Proposition 3, \cite{K1}). The case $z=0$ does not follow the proof in \cite{K1}. 

If $z=0$, we prove directly that $v(0) > u(0)$. Let us suppose for simplicity that $u$ and $v$ are continuous functions. Then for any small $\varepsilon >0$ we can find the probe functions 
$\phi_\varepsilon^{(1)}, \phi_\varepsilon^{(2)} \in \mathbf{C}_1[0,1]$ such that
$$
\begin{aligned}
x_0^{(1)}&:=\mbox{argmax}(u(x) - \phi_\varepsilon^{(1)}(x)) < \varepsilon \\ x_0^{(2)}&:=\mbox{argmin}(v(x) - \phi_\varepsilon^{(2)}(x)) < \varepsilon
\end{aligned}
$$
and 
$$
\begin{aligned}
&u(x_0^{(1)}) - \beta H \left(x_0^{(1)}, \frac{\ed}{\ed x} \phi_\varepsilon^{(1)} (x_0^{(1)})\right) - h(x_0^{(1)}) \le 0, \\
&v(x_0^{(2)}) - \beta H\left(x_0^{(2)}, \frac{\ed}{\ed x} \phi_\varepsilon^{(2)} (x_0^{(2)})\right) - h(x_0^{(2)}) \ge 0.
\end{aligned}
$$
Getting $\varepsilon \to 0$ we obtain that $v(0)>u(0)$. The case of semicontinuous $u$ and $v$ is proved in a similar way adding some conditions for probe functions. It finishes the proof of the proposition. $\Box$

\section*{Appendix 2: control theory}\label{A4}

We need to verify if this semigroup is well-defined: we need to check Conditions 8.9 and 8.11 from \cite{Kurtz2}. To verify these conditions the following definition of a control equation is required. Based on trivial relation (deterministic control)
$$
f(\Ab{\gamma}(t)) - f(\Ab{\gamma}(0)) = \int_0^t f^\prime(\Ab{\gamma}(s)) \dot{\Ab{\gamma}}(s) ds,
$$
the control theory introduce the notion of \textit{relaxed control equation}, where the deterministic velocity $\dot{\Ab{\gamma}}$ is substituted by a (control) measure $\lambdabar$, see definition below.

\begin{definition} (\cite{Kurtz2}, Definition 8.1, or \cite{Ber}, Definition 3.5) Let $U$ and $E$ be complete and separable metric spaces. Let $A:\ Dom(A) \subset B(E) \to  M(E \times U)$ be a single valued linear operator.
	Let $\mathcal{M}_m(U)$ be the space of Borel measures $\lambdabar$ on $U\times[0, 1]$ satisfying $\lambdabar (U \times [0, t]) = t$ for all $t\in [0, 1]$. The measure $\lambdabar$ is known as a {\rm relaxed control}. We say that the pair $(\mathbf{x}, \lambdabar) \in D_E[0, 1]\times \mathcal{M}_m(U)$ satisfies
	the {\rm relaxed control equation} for $A$ if and only if:
	\begin{enumerate}
	\item for any $f\in Dom(A)$ any $t\in [0,1]$
	$$
	\iint_{U\times[0,1]} \left| A(f)(\mathbf{x}(s), u)\right| \lambdabar(du \times ds) < \infty,
	$$
	\item for any $f\in Dom(A)$ any $t\in [0,1]$
	$$
	f(\mathbf{x}(t)) - f(\mathbf{x}(0)) = \iint_{U\times[0,1]} \left| A(f)(\mathbf{x}(s), u)\right| \lambdabar(du \times ds).
	$$
	\end{enumerate}
\end{definition} 
Observe that the deterministic control equation is a relaxed control equation assuming
\begin{equation}\label{det.control}
A(f)(x, v) = f^\prime(x)v \mbox{ and }\lambdabar(dv,ds) = \delta_{v(s)}(dv) ds.
\end{equation}
Moreover, in our case the operator $H$ is $H(x,f^\prime(x))$ for each $x\in [0,1]$, and \reff{Hxp}
 can be rewritten in terms of operator $A$ as
\begin{equation}\label{Hxp1}
	H(f)(x) = H(x,f^\prime(x)) = \sup_{u} \left( A(f)(x,u) - L(x,u) \right).
\end{equation}
  
Let $\mathcal Y$ be the set of all such pairs $(\mathbf{x}, \lambdabar)$. Instead of \reff{Nisio1} the relaxed controls determine the Nisio semigroup as
\begin{equation}\label{Nisio2}
	\mathbf{V}(t)f(x) = \sup_{(\mathbf{x},\lambdabar) \in \mathcal{Y},\ \mathbf{x}(0)=x} f(\mathbf{x}(t)) - \iint_{U\times[0,1]} L(\mathbf{x}(s),u) \lambdabar(du \times ds)
\end{equation}

Theorem~3.1 (section 4.3.2, see \cite{Ber}) provides the variational and control representation of the rate function. 

We finish the section with Conditions 8.9 and 8.11 from \cite{Kurtz2}.

\vspace{0.5cm}
\noindent
\textbf{Condition 8.9 \cite{Kurtz2}} 
\begin{enumerate}
	\item[(1)] $A\subset C_b(E)\times C(E\times U)$ is single-valued and $Dom(A)$ separates points.
	\item[(2)] $\Gamma \subset E\times U$ is closed, and for each $x_0\in E$, there exists $(\mathbf{x},\lambdabar)$ such that $\mathbf{x}(0)=x_0$ and $$\int_{U\times [0,t]} \mathbb{I}_\Gamma(\mathbf{x}(s),u) \lambdabar(du\times ds) =t,$$ where $\mathbb{I}_B(\cdot)$ is an indicator of a set $B$.
	\item[(3)] $L(x,u): E\times U \to [0,\infty]$ is a lower semicontinuous function, and for each $c\in[0,\infty)$ and compact $K\subset E$,
	$$
	\{ (x,u) \in \Gamma: \ L(x,u) \le c\} \cap (K\times U)
	$$
	is relatevely compact.  
	\item[(4)] For each compact $K\subset E, T > 0$, and $0\le M <\infty$, there exists a compact
	$\hat{K}\equiv \hat{K}(K,T,M) \subset E$ such that $\mathbf{x}(0) \in K$, 
	$$
	\int_{U\times [0,t]} \mathbb{I}_\Gamma(\mathbf{x}(s),u) \lambdabar(du\times ds) =t,
	\ \mbox{ and }\ 
	\int_{U\times [0,T]} L(\mathbf{x}(s), u)\lambdabar(du\times ds) \le M
	$$
	imply $\mathbf{x}(t)\in \hat{K}, \ 0\le t \le T$.
	\item[(5)] For each $f\in Dom(A)$ and compact set $K\subset E$, there exists a right continuous, nondecreasing function $\psi_{f,K}:\ [0,\infty) \to [0,\infty)$ such that
	$$
	|Af(x,u)| \le \psi_{f,K}(L(x,u)), \ \ \forall\ (x,u)\in \Gamma \cap (K\times U),
	$$
	and 
	$$
	\lim_{r\to\infty} \frac{\psi_{f,K}(r)}{r} =0.
	$$
\end{enumerate}

\vspace{0.5cm}

\noindent
\textbf{Condition 8.11 \cite{Kurtz2}} For each $x_0\in E$ and each $f\in Dom(\mathcal H)$, there exists $(\mathbf{x},\lambdabar)$ with $\mathbf{x}(0)=x_0$ and $$\int_{U\times [0,t]} \mathbb{I}_\Gamma(\mathbf{x}(s),u) \lambdabar(du\times ds) =t,$$ such that for $0\le t_1 < t_2$
$$
\int_{t_1}^{t_2} \mathcal{H}f(\mathbf{x}(s))ds \le \int_{U\times (t_1,t_2]} \left( 
Af(\mathbf{x}(s),u) - L(\mathbf{x}(s),u) \right) \lambdabar(ds\times du).
$$
These conditions are easily verified in our case. See also \cite{Ber}.

\section*{Appendix 3}\label{A1}

Consider some cases.

\vspace{0.5cm}
\noindent
{\bf Case $\Ab{\gamma}(0)=\Ab{\gamma}(T)=y$.} In this case the solution is constant, $\Ab{\gamma}(t) \equiv \Ab{\gamma}(0)$.

\begin{figure}[h]
	\caption{Here we plot the optimal paths for the case $T=2$ and boundaries conditions $\Ab{\gamma}(0)=0$ and $\Ab{\gamma}(0)< \Ab{\gamma}(2) \in \{01, 0.3, 0.5, 0.9\}$. The upper small plot shows the whole parabolas -- solutions of system of Hamiltonian equations. }	
	\centering 
	\includegraphics[width=15cm]{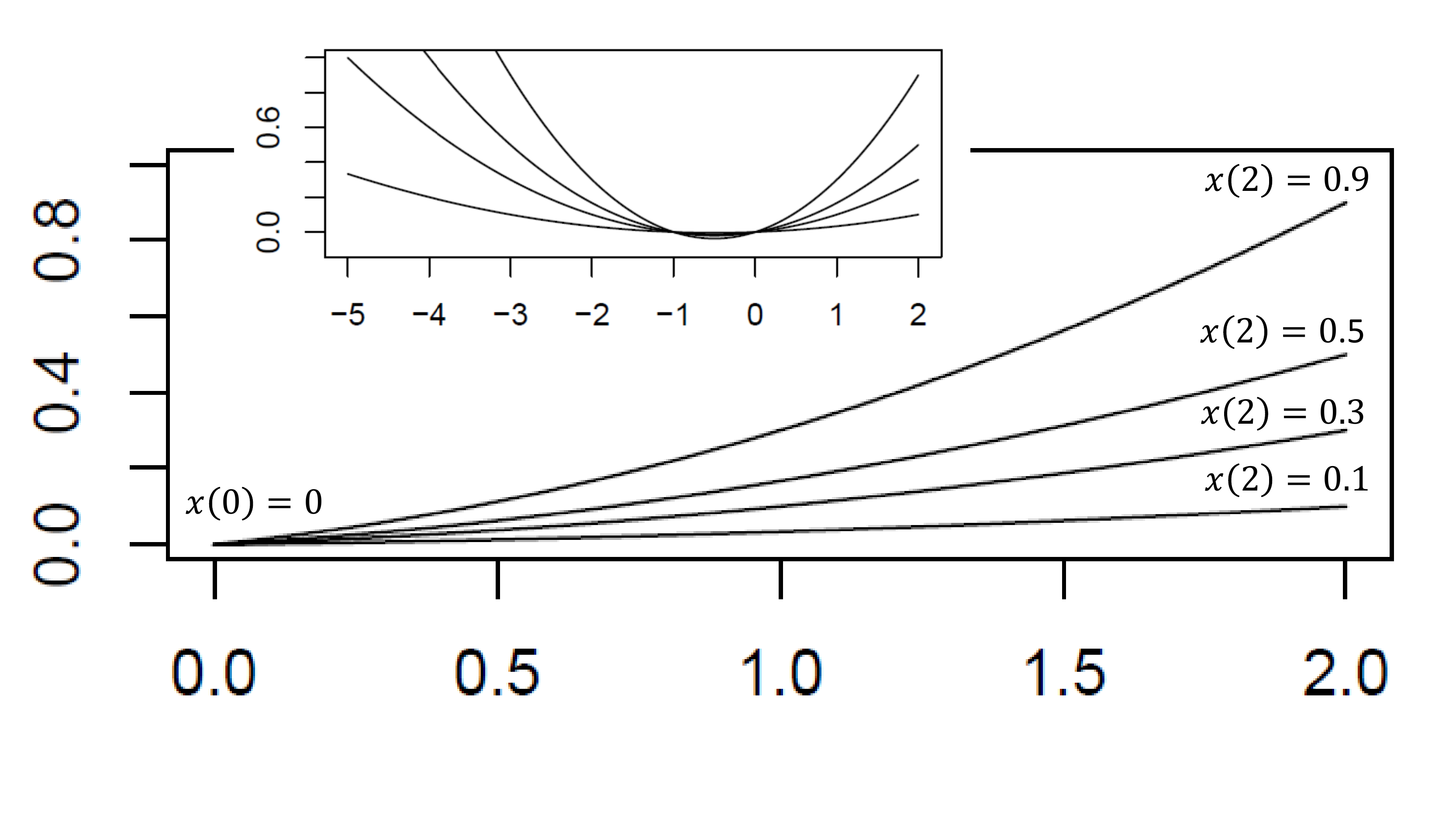} 
	\label{fig2}
\end{figure}

\vspace{0.5cm}
\noindent
{\bf Case $0 \le \Ab{\gamma}(0) < \Ab{\gamma}(T)$.} Consider first the case when $\Ab{\gamma}(0)=0$ and $\Ab{\gamma}(T)\ne 0$. In this case the formulas are simple
$$
\begin{aligned}
	&c_1=0, \ \ c_2 = \frac{\Ab{\gamma}(T)}{\lambda T(\lambda T+1)} \\
	&x(t) = \lambda t (\lambda t + 1) \frac{\Ab{\gamma}(T)}{\lambda T(\lambda T+1)}
\end{aligned}
$$
In Figure~\ref{fig2} we plot four the function $\Ab{\gamma}(t)$ for case when $T=2$ and $\Ab{\gamma}(2)=0.1, \Ab{\gamma}(2)=0.3, \Ab{\gamma}(2)=0.5$ and $\Ab{\gamma}(2)=0.9$.

\begin{figure}[h]
	\caption{\it Here we plot the optimal paths for the case $T=2$ and boundaries conditions $\Ab{\gamma}(0)=0$ and $\Ab{\gamma}(2) \in \{0.0, 0.3, 0.5, 0.7, 1.0\}$. In order to show whole parabolas, we add the upper small plot, which shows the same parabolas in larger "time" interval $t\in[-10,10]$. Rectangle outlined in red is the area of optimal trajectories $\Ab{\gamma}(t)$: $[0,1]\times [0,T]$.}	
	\centering 
	\includegraphics[width=15cm]{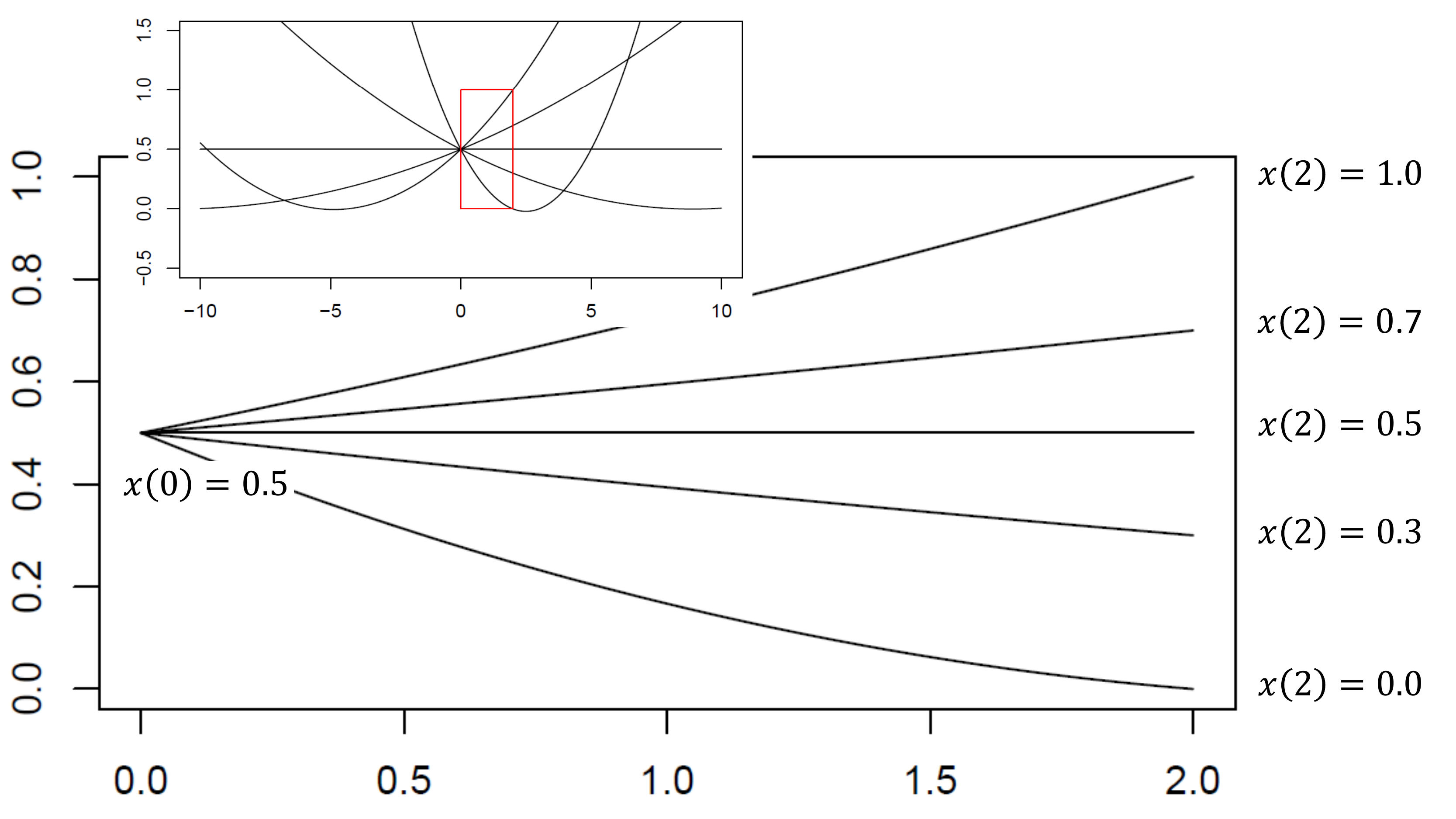} 
	\label{fig3}
\end{figure}

We should also note that when values $\Ab{\gamma}(0)$ and $\Ab{\gamma}(T)$ approach each other, then the parabola $\Ab{\gamma}(t)$ transforms into horizontal line $\Ab{\gamma}(t) \equiv \Ab{\gamma}(0)$. 

If $\Ab{\gamma}(0)>0$ and $\Ab{\gamma}(0)< \Ab{\gamma}(T)$, then $c_1$ we find as the solution of quadratic equation
\begin{equation}\label{qua-eq}
	c_1^2 + c_1 \left( \frac{2\lambda T\Ab{\gamma}(0)}{\Ab{\gamma}(T)-\Ab{\gamma}(0)} -1 \right) - \frac{\lambda T (\lambda T +1)\Ab{\gamma}(0)}{\Ab{\gamma}(T)-\Ab{\gamma}(0)} = 0, \ \ c_2 = \frac{\Ab{\gamma}(0)}{c_1(c_1-1)}
\end{equation}
The value of $c_1$ which provides the admissible function $\Ab{\gamma}(t)$, i.e. $\Ab{\gamma}(t) \in [0,1]$ for all $t\in [0,T]$ is
$$
c_1 = \frac{1}{2} - \frac{\lambda T\Ab{\gamma}(0)}{\Ab{\gamma}(T)-\Ab{\gamma}(0)} -
\sqrt{\left( \frac{1}{2} - \frac{\lambda T\Ab{\gamma}(0)}{x(T)-\Ab{\gamma}(0)} \right)^2 + \frac{\lambda T (\lambda T + 1)\Ab{\gamma}(0)}{\Ab{\gamma}(T)-\Ab{\gamma}(0)}}, \ \ c_2 = \frac{\Ab{\gamma}(0)}{c_1(c_1-1)}
$$  
For the time moment $T=2$ see Figure~\ref{fig3}. The here considered case corrsponds to the boundary conditions $x(0)=0.5, x(2)=0.7$ and $x(0)=0.5, x(2)=1$ in the plot.

\vspace{0.5cm}
\noindent
{\bf Case $0 \le \Ab{\gamma}(T) < \Ab{\gamma}(0)$.} Consider first the case when $\Ab{\gamma}(T)=0$ and $\Ab{\gamma}(0)\ne 0$. In this case constants which defines the admissible trajectories are 
$$
c_1=\lambda T +1, \ \ c_2 = \frac{\Ab{\gamma}(0)}{\lambda T(\lambda T+1)} \mbox{ and }
\Ab{\gamma}(t) = \lambda(T-t)(\lambda(T-t) + 1) \frac{x(0)}{\lambda T(\lambda T+1)}
$$
Example can be found in Figure~\ref{fig3} for the case $T=2$ and boundary condition $\Ab{\gamma}(0)=1/2$ and $\Ab{\gamma}(2)=0$. The optimal trajectories are parabolas again, and they are symmetric version of the Figure~\ref{fig2}. 

If $\Ab{\gamma}(T)>0$ and $\Ab{\gamma}(0) > \Ab{\gamma}(T)$, then $c_1$ we find as the solution of the same quadratic equation \reff{qua-eq}.
The root of \reff{qua-eq} which provides the admissible function $\Ab{\gamma}(t)$ is
$$
c_1 = \frac{1}{2} - \frac{\lambda Tx(0)}{\Ab{\gamma}(T)-\Ab{\gamma}(0)} +
\sqrt{\left( \frac{1}{2} - \frac{\lambda T\Ab{\gamma}(0)}{\Ab{\gamma}(T)-\Ab{\gamma}(0)} \right)^2 + \frac{\lambda T (\lambda T + 1)\Ab{\gamma}(0)}{\Ab{\gamma}(T)-\Ab{\gamma}(0)}}, \ \ c_2 = \frac{\Ab{\gamma}(0)}{c_1(c_1-1)}
$$  
As an example, see Figure:\ref{fig3} with $T=2$, $\Ab{\gamma}(0)=0.5, \Ab{\gamma}(2)=0.3$ and $\Ab{\gamma}(0)=0.5, \Ab{\gamma}(2)=0.0$. 

Observe that extremal points of the corresponding parabolas are greater then $T$, when these extremal points in the case $\Ab{\ga}(0)<\Ab{\ga}(T)$ are negative.


\vspace{5cm}

\end{document}